\documentclass[12pt]{article}
\usepackage{amsfonts}
\usepackage{full page,
amssymb, amscd, amsmath,graphicx, 
}

 \title{{\bf Twisted representations of vertex operator algebras associated to affine Lie algebras}}
 \author{Jinwei Yang}
    \date{}
    \begin{document}
    \bibliographystyle{alpha}
    \maketitle
\newtheorem{thm}{Theorem}[section]
\newtheorem{defn}[thm]{Definition}
\newtheorem{prop}[thm]{Proposition}
\newtheorem{cor}[thm]{Corollary}
\newtheorem{lemma}[thm]{Lemma}
\newtheorem{rema}[thm]{Remark}
\newtheorem{app}[thm]{Application}
\newtheorem{prob}[thm]{Problem}
\newtheorem{conv}[thm]{Convention}
\newtheorem{conj}[thm]{Conjecture}
\newtheorem{cond}[thm]{Condition}
    \newtheorem{exam}[thm]{Example}
\newtheorem{assum}[thm]{Assumption}
     \newtheorem{nota}[thm]{Notation}
\newcommand{\halmos}{\rule{1ex}{1.4ex}}
\newcommand{\pfbox}{\hspace*{\fill}\mbox{$\halmos$}}

\newcommand{\nn}{\nonumber \\}

 \newcommand{\res}{\mbox{\rm Res}}
 \newcommand{\ord}{\mbox{\rm ord}}
\renewcommand{\hom}{\mbox{\rm Hom}}
\newcommand{\edo}{\mbox{\rm End}\ }
 \newcommand{\pf}{{\it Proof.}\hspace{2ex}}
 \newcommand{\epf}{\hspace*{\fill}\mbox{$\halmos$}}
 \newcommand{\epfv}{\hspace*{\fill}\mbox{$\halmos$}\vspace{1em}}
 \newcommand{\epfe}{\hspace{2em}\halmos}
\newcommand{\nord}{\mbox{\scriptsize ${\circ\atop\circ}$}}
\newcommand{\wt}{\mbox{\rm wt}\ }
\newcommand{\swt}{\mbox{\rm {\scriptsize wt}}\ }
\newcommand{\lwt}{\mbox{\rm wt}^{L}\;}
\newcommand{\rwt}{\mbox{\rm wt}^{R}\;}
\newcommand{\slwt}{\mbox{\rm {\scriptsize wt}}^{L}\,}
\newcommand{\srwt}{\mbox{\rm {\scriptsize wt}}^{R}\,}
\newcommand{\clr}{\mbox{\rm clr}\ }
\newcommand{\tr}{\mbox{\rm Tr}}
\newcommand{\C}{\mathbb{C}}
\newcommand{\Z}{\mathbb{Z}}
\newcommand{\R}{\mathbb{R}}
\newcommand{\Q}{\mathbb{Q}}
\newcommand{\N}{\mathbb{N}}
\newcommand{\CN}{\mathcal{N}}
\newcommand{\F}{\mathcal{F}}
\newcommand{\I}{\mathcal{I}}
\newcommand{\V}{\mathcal{V}}
\newcommand{\one}{\mathbf{1}}
\newcommand{\BY}{\mathbb{Y}}
\newcommand{\ds}{\displaystyle}

        \newcommand{\ba}{\begin{array}}
        \newcommand{\ea}{\end{array}}
        \newcommand{\be}{\begin{equation}}
        \newcommand{\ee}{\end{equation}}
        \newcommand{\bea}{\begin{eqnarray}}
        \newcommand{\eea}{\end{eqnarray}}
         \newcommand{\lbar}{\bigg\vert}
        \newcommand{\p}{\partial}
        \newcommand{\dps}{\displaystyle}
        \newcommand{\bra}{\langle}
        \newcommand{\ket}{\rangle}

        \newcommand{\ob}{{\rm ob}\,}
        \renewcommand{\hom}{{\rm Hom}}

\newcommand{\A}{\mathcal{A}}
\newcommand{\Y}{\mathcal{Y}}

\newcommand{\dlt}[3]{#1 ^{-1}\delta \bigg( \frac{#2 #3 }{#1 }\bigg) }

\newcommand{\dlti}[3]{#1 \delta \bigg( \frac{#2 #3 }{#1 ^{-1}}\bigg) }

 \makeatletter
\newlength{\@pxlwd} \newlength{\@rulewd} \newlength{\@pxlht}
\catcode`.=\active \catcode`B=\active \catcode`:=\active
\catcode`|=\active
\def\sprite#1(#2,#3)[#4,#5]{
   \edef\@sprbox{\expandafter\@cdr\string#1\@nil @box}
   \expandafter\newsavebox\csname\@sprbox\endcsname
   \edef#1{\expandafter\usebox\csname\@sprbox\endcsname}
   \expandafter\setbox\csname\@sprbox\endcsname =\hbox\bgroup
   \vbox\bgroup
  \catcode`.=\active\catcode`B=\active\catcode`:=\active\catcode`|=\active
      \@pxlwd=#4 \divide\@pxlwd by #3 \@rulewd=\@pxlwd
      \@pxlht=#5 \divide\@pxlht by #2
      \def .{\hskip \@pxlwd \ignorespaces}
      \def B{\@ifnextchar B{\advance\@rulewd by \@pxlwd}{\vrule
         height \@pxlht width \@rulewd depth 0 pt \@rulewd=\@pxlwd}}
      \def :{\hbox\bgroup\vrule height \@pxlht width 0pt depth
0pt\ignorespaces}
      \def |{\vrule height \@pxlht width 0pt depth 0pt\egroup
         \prevdepth= -1000 pt}
   }
\def\endsprite{\egroup\egroup}
\catcode`.=12 \catcode`B=11 \catcode`:=12 \catcode`|=12\relax
\makeatother

\def\hboxtr{\FormOfHboxtr} 
\sprite{\FormOfHboxtr}(25,25)[0.5 em, 1.2 ex] 

:BBBBBBBBBBBBBBBBBBBBBBBBB | :BB......................B |
:B.B.....................B | :B..B....................B |
:B...B...................B | :B....B..................B |
:B.....B.................B | :B......B................B |
:B.......B...............B | :B........B..............B |
:B.........B.............B | :B..........B............B |
:B...........B...........B | :B............B..........B |
:B.............B.........B | :B..............B........B |
:B...............B.......B | :B................B......B |
:B.................B.....B | :B..................B....B |
:B...................B...B | :B....................B..B |
:B.....................B.B | :B......................BB |
:BBBBBBBBBBBBBBBBBBBBBBBBB |

\endsprite
\def\shboxtr{\FormOfShboxtr} 
\sprite{\FormOfShboxtr}(25,25)[0.3 em, 0.72 ex] 

:BBBBBBBBBBBBBBBBBBBBBBBBB | :BB......................B |
:B.B.....................B | :B..B....................B |
:B...B...................B | :B....B..................B |
:B.....B.................B | :B......B................B |
:B.......B...............B | :B........B..............B |
:B.........B.............B | :B..........B............B |
:B...........B...........B | :B............B..........B |
:B.............B.........B | :B..............B........B |
:B...............B.......B | :B................B......B |
:B.................B.....B | :B..................B....B |
:B...................B...B | :B....................B..B |
:B.....................B.B | :B......................BB |
:BBBBBBBBBBBBBBBBBBBBBBBBB |

\endsprite

\vspace{2em}

\begin{abstract}
In this paper, we prove the categories of lower bounded twisted modules of positive integer levels for simple vertex operator algebras associated with affine Lie algebras and general automorphisms are semisimple, using the twisted generalization of Zhu's algebra for these vertex operator algebras, constructed in \cite{HY}. We also show that the category of lower bounded twisted modules for a general automorphism is equivalent to the category of lower bounded twisted modules for the corresponding diagram automorphism.
\end{abstract}

\renewcommand{\theequation}{\thesection.\arabic{equation}}
\renewcommand{\thethm}{\thesection.\arabic{thm}}
\setcounter{equation}{0} \setcounter{thm}{0} \date{}
\maketitle

\section{Introduction}
Orbifold conformal field theories play an important role in mathematics and physics. They
are examples of conformal field theories constructed from known conformal field theories and
automorphisms of these known ones.
Mathematically, the study of orbifold conformal field theories
can be reduced to the study of the twisted
representation theory of vertex operator algebras, that is,
the theory of representations of vertex operator algebras twisted by automorphisms.
In fact, the first example of orbifold conformal field theories is the conformal field
theory corresponding to the moonshine module vertex operator algebra constructed by
Frenkel, Lepowsky and Meurman \cite{FLM}.
Orbifold theory was then studied in the physics literature (including \cite{DHVW}, \cite{DVVV}) and as twisted module theory for
a vertex operator algebra in the mathematics literature (see for
example \cite{DL}, \cite{DLM}, \cite{L}, \cite{H}, \cite{B}, \cite{HY}).

Automorphisms of a vertex operator algebra can be of finite or infinite orders.
An element of the Monster group gives an automorphism of the moonshine module
vertex operator algebra \cite{FLM}. Such an automorphism
is of finite order. However, for most of the classes of vertex operator algebras, the automorphism groups are groups of Lie type. The full automorphism groups of Heisenberg vertex operator algebras, affine vertex operator algebras, symplectic fermion algebras, lattice vertex operator algebras, triplet vertex operator algebras preserving a natural conformal structure are all groups of Lie type (see \cite{FZ}, \cite{FLM}, \cite{Lin}, \cite{CLin}, \cite{DN}, \cite{ALM}). An element in these full automorphism groups is in general of infinite order. Moreover, in general it might not act on the vertex operator algebra semisimply.

In \cite{H}, inspired by the logarithmic conformal field theory, Y.-Z. Huang generalized
the notion of (lower bounded) twisted module for a vertex operator algebra with a
finite-order automorphism to the notion of (lower bounded) twisted module
for a vertex operator algebra with a general automorphism not
necessarily of finite order.  One important new feature of twisted modules in the general case is that twisted vertex operators might involve the logarithm of the variable. We call the twisted modules in the general case the {\it logarithmic} twisted modules.

In the present paper, we study the logarithmic twisted modules for the vertex operator algebras associated to the affine Lie algebra of a finite dimensional Lie algebra. The automorphisms of these classes of vertex operator algebras arise from the automorphisms of the corresponding finite dimensional Lie algebras, which, in the case of a compact simple Lie algebra, are products of inner automorphisms and diagram automorphisms (see for example \cite{He}). For finite order such automorphisms, Haisheng Li  (\cite{L}) gave an equivalent condition between the irreducible lower bounded twisted modules for the simple affine vertex operator algebras and the standard modules for the twisted affine Lie algebras. He also determined the complete reducibility for these lower bounded twisted modules.

In general, for an arbitrary automorphism $g$ of a vertex operator algebra $V$, jointly with Y.-Z. Huang, the author constructed an associative algebra $A_g(V)$, called {\it $g$-twisted Zhu's algebra}, together with functors between the categories of suitable logarithmic $g$-twisted $V$-modules and $A_g(V)$-modules, generalizing the constructions by Dong, Li and Mason in \cite{DLM} (see also \cite{KW} for a definition of associative algebras for vertex operator superalgebras). Using these constructions, we can reduce the study of lower bounded $g$-twisted modules to the study of the associative algebra modules. To determine the associative algebra $A_g(V)$ is the first key step in this process, we will describe the twisted Zhu's algebras for the vertex operator algebras associated to the affine Lie algebras in this paper.

For a finite dimensional Lie algebra $\mathfrak{g}$ and a positive integer $\ell$. Let $V_\mathfrak{g}(0, \ell)$ be the universal affine vertex operator algebra of level $\ell$, and let $L_\mathfrak{g}(0, \ell)$ be the simple quotient of $V_{\mathfrak{g}}(0, \ell)$. Let $\varphi$ be an automorphism of $\mathfrak{g}$. We first describe the $\varphi$-twisted Zhu's algebra for $V_{\mathfrak{g}}(0, \ell)$ and for $L_\mathfrak{g}(0, \ell)$, then use the bijections on the isomorphism classes of irreducible lower bounded $\varphi$-twisted modules and of irreducible twisted Zhu's algebra modules, we give a necessary condition for irreducible lower bounded twisted modules. The complete reducibility of lower bounded logarithmic twisted $L_\mathfrak{g}(0, \ell)$-modules is also a consequence of the description of the twisted Zhu's algebras.

An important feature of the twisted Zhu's algebras for $V_\mathfrak{g}(0, \ell)$ is that they are independent of the inner automorphisms (see also \cite{B}). For the simple vertex operator algebra $L_\mathfrak{g}(0, \ell)$, we first show that for the automorphisms induced from inner automorphisms of $\mathfrak{g}$, the twisted Zhu's algebras for $L_\mathfrak{g}(0, \ell)$ are all the same as the usual (untwisted) Zhu's algebra. As a consequence, we determine the list of irreducible lower bounded twisted $L_\mathfrak{g}(0, \ell)$-modules for the inner automorphisms. In general, let $\mu$ be the diagram automorphism involved in the automorphism $\varphi$ of $\mathfrak{g}$. We prove that the category of lower bounded $\varphi$-twisted modules and the category of lower bounded $\mu$-twisted modules are equivalent. As a consequence, the twisted Zhu's algebras for $L_\mathfrak{g}(0, \ell)$ are independent of the inner automorphisms.

The present paper is organized as follows: In Section 2, we recall the definitions and properties of twisted logarithmic modules. In Section 3, we recall the constructions of twisted Zhu's algebra and functors between twisted logarithmic modules and twisted Zhu's algebra modules. We recall the notions and properties of vertex operator algebras associated to the affine Lie algebras in Section 4. In Section 5 and 6, we prove the main results for general automorphisms.

\paragraph{Acknowledgments}
We would like to thank Bojko Bakalov, Katrina Barron, Yi-Zhi Huang, Haisheng Li and Robert McRae for valuable discussions.

\section{Logarithmic twisted modules}
Let $(V, Y, {\bf 1}, \omega)$ be a vertex
operator algebra and let $g$ be an automorphism of $V$. In this section, we recall some definitions and properties of logarithmic $g$-twisted $V$-modules from \cite{HY} (see also \cite{H} and \cite{B}). We shall use the convention that $\log z=\log |z|+\arg z \sqrt{-1}$
where $0\le \arg z<2\pi$. We shall also use $l_{p}(z)$ to denote $\log z+2 \pi p\sqrt{-1}$ for $p\in \Z$.
Besides the usual notations $\C$, $\R$, $\Z$, $\Z_{+}$ and $\N$ for the sets of complex numbers,
real numbers, integers, positive integers and natural numbers, respectively, we also use $\mathbb{I}$ and
$\overline{\C}_{+}$ to denote the set of imaginary numbers and the closed right half complex plane.

For $\alpha \in \C/\Z$, we denote by
$$
V^{[\alpha]} = \{u \in V\;|\; (g-e^{2\pi\sqrt{-1}\alpha})^{\Lambda} u =0 \;\text{for some}\; \Lambda\in \Z_{+}\}
$$
the generalized eigenspace with eigenvalue $e^{2\pi \sqrt{-1}\alpha}$ for $g$.
For $\alpha\in \C/\Z$,
there is a unique $a\in [0, 1)+\mathbb{I}$ such that $a+\Z=\alpha$.
We call $a\in [0, 1)+\mathbb{I}$ a
{\it $g$-weight of $V$} if $V^{[a+\Z]} \neq 0$.
We use $P(V)$ to denote the
set of all the $g$-weights of $V$.  Then
$V=\coprod_{n\in \C, a\in P(V)}V_{(n)}^{[a+\Z]}$.

\begin{defn}\label{d}
{\rm  A {\it  $\overline{\C}_{+}$-graded weak $g$-twisted
$V$-module} is a $\overline{\C}_{+}\times \mathbb{C}/\mathbb{Z}$-graded
vector space $W = \coprod_{n \in \overline{\C}_{+}, \alpha \in
\mathbb{C}/\mathbb{Z}} W_{n}^{[\alpha]}$ (graded by {\it $\overline{\C}_{+}$-degrees} and {\it $g$-weights})
equipped with a linear
map
\begin{align*}
Y^g: &V\otimes W \rightarrow W\{x\}[{\rm log} x],\\
&v \otimes w \mapsto Y^g(v, x)w
\end{align*}
 and an action of $g$ satisfying the following conditions:
\begin{enumerate}

\item The {\it equivariance property}: For $p \in \mathbb{Z}$, $z
\in \mathbb{C}^{\times}$, $v \in V$ and $w \in W$, $Y^{g; p + 1}(gv,
z)w = Y^{g; p}(v, z)w$, where for $p \in \mathbb{Z}$,
$$Y^{g; p}(v, z)w=Y^{g}(v, x)w\lbar_{x^{n}=e^{nl_{p}(z)},\; \log x=l_{p}(z)}$$
is the $p$-th analytic branch of $Y^g$.

\item The {\it identity property}: For $w \in W$, $Y^g({\bf 1}, x)w
= w$.

\item  The {\it duality property}: Let $W^{'} = \coprod_{n \in \overline{\C}_{+}, \alpha \in
\mathbb{C}/\mathbb{Z}} \left(W_{n}^{[\alpha]}\right)^{\ast}$ and, for $n \in
\mathbb{C}$, $\pi_n : W \rightarrow W_{n}=\coprod_{\alpha\in \C/\Z}W_{n}^{[\alpha]}$ be the projection. For
any $u, v \in V$, $w \in W$ and $w^{'} \in W^{'}$, there exists a
multivalued analytic function of the form
\[
f(z_1, z_2) = \sum_{i,
j, k, l = 0}^N a_{ijkl}z_1^{m_i}z_2^{n_j}({\rm log}z_1)^k({\rm
log}z_2)^l(z_1 - z_2)^{-t}
\]
for $N \in \mathbb{N}$, $m_1, \dots,
m_N$, $n_1, \dots, n_N \in \mathbb{C}$ and $t \in \mathbb{Z}_{+}$,
such that the series
\[
\langle w^{'}, Y^{g; p}(u, z_1)Y^{g;p}(v,
z_2)w\rangle = \sum_{n \in \mathbb{C}}\langle w^{'}, Y^{g; p}(u,
z_1)\pi_nY^{g; p}(v, z_2)w\rangle,
\]
\[
\langle w^{'}, Y^{g; p}(v,
z_2)Y^{g;p}(u, z_1)w\rangle = \sum_{n \in \mathbb{C}}\langle w^{'},
Y^{g; p}(v, z_2)\pi_nY^{g; p}(u, z_1)w\rangle,
\]
\[
\langle w^{'}, Y^{g; p}(Y(u, z_1 - z_2)v,
z_2)w\rangle = \sum_{n \in \mathbb{C}}\langle w^{'}, Y^{g;
p}(\pi_nY(u, z_1 - z_2)v, z_2)w\rangle
\]
are absolutely convergent in
the regions $|z_1| > |z_2| > 0$, $|z_2| > |z_1| > 0$, $|z_2| > |z_1
- z_2| > 0$, respectively, and are convergent to the branch
\[
 \sum_{i, j, k, l = 0}^N
a_{ijkl}e^{m_il_p(z_1)}e^{n_jl_p(z_2)}l_p(z_1)^kl_p(z_2)^l(z_1 -
z_2)^{-t}
\]
of $f(z_1, z_2)$ when $\arg z_{1}$ and $\arg z_{2}$ are sufficiently close (more precisely,
when $|\arg z_{1}-\arg z_{2}|<\frac{\pi}{2}$).

\item The {\it $\overline{\C}_{+}$- and $g$-grading conditions}: For $v\in V_{(m)}$
 and $w\in W_{p}=\coprod_{\alpha\in \C/\Z}W_{p}^{[\alpha]}$ where
$m\in \Z$ and $p\in \overline{\C}_{+}$, write
$Y^{g}(v, x)w=\sum_{k=0}^{N}\sum_{n\in \C}Y^{g}_{n, k}(v)wx^{-n-1}(\log x)^{k}.$
Then $Y^{g}_{n, k}(v)w$ is $0$ when $m-n-1+p\not\in \overline{\C}_{+}$ (or $\Re(m-n-1+p)<0$),
is  in $W_{m-n-1+p}$ when $m-n-1+p\in \overline{\C}_{+}$ and for $r\in \R$,
$\coprod_{n\in \mathbb{I}}W_{n}$ is equal to the subspace of $W$ consisting
of $w\in W$ such that when  $Y^{g}_{n, k}(v)w=0$ for
$v\in V_{(m)}$, $n\in \C$, $k=0, \dots, N$, $m-n-1\not\in \overline{\C}_{+}$.
For $\alpha \in \mathbb{C}/\mathbb{Z}$, $w \in
W^{[\alpha]}=\coprod_{n\in \overline{\C}_{+}}W_{n}^{[\alpha]}$, there exists $\Lambda \in \mathbb{Z}_{+}$
such that $ (g -
e^{2\pi\sqrt{-1}\alpha})^{\Lambda}w = 0$. Moreover, $gY^{g}(u, x)v=Y^{g}(gu, x)gv$.

\item The $L(-1)$-{\it derivative property}: For $v \in V$,
\[
\frac{d}{dx}Y^g(v, x) = Y^g(L(-1)v, x).
\]
\end{enumerate}

A {\it lower bounded logarithmic $g$-twisted $V$-module} or simply a {\it lower bounded $g$-twisted $V$-module}
is a $\overline{\C}_{+}$-graded weak $g$-twisted
$V$-module $W$ together with a decomposition of $W$ as a direct sum $W=\coprod_{n\in \C}W_{[n]}$
of generalized eigenspaces $W_{[n]}$ with eigenvalues $n\in \C$ for the operator
$L^{g}(0)=\res_{x}xY^g(\omega, x) $ such that
for each $n\in \C$ and each $\alpha\in \C/\Z$, $W^{[\alpha]}_{[n + l]} =W_{[n+l]}\cap W^{[\alpha]}= 0$ for
sufficiently negative real number $l$. A lower bounded logarithmic $g$-twisted $V$-module $W$
is said to be {\it strongly $\C/\Z$-graded} or {\it grading-restricted} if it is lower
bounded and for each $n \in
\mathbb{C}$, $\alpha \in \mathbb{C}/\mathbb{Z}$,
$\dim W^{[\alpha]}_{[n]} =\dim W_{[n]}\cap W^{[\alpha]}< \infty$.}
\end{defn}

Let $W=\coprod_{n\in \overline{\C}_{+},\alpha\in \C/\Z}W_{n}^{[\alpha]}$ be a
$\overline{\C}_{+}$-graded weak $g$-twisted
$V$-module with the twisted vertex operator map $Y^{g}$.
For $\alpha\in \C/\Z$, let $W^{[\alpha]}=\coprod_{n\in \overline{\C}_{+}}W_{n}^{[\alpha]}$.
When $W$ is a (lower bounded) $g$-twisted $V$-module,  we also have $W^{[\alpha]}=
\coprod_{n\in \C}W_{[n]}^{[\alpha]}$.
As in the case of $V$, we call $a\in [0, 1)+\mathbb{I}$ a
{\it $g$-weight of $W$} if $W^{[a+\Z]} \neq 0$.
Let $P(W)$ be the set of all $g$-weights of $W$. Then
$W=\coprod_{n\in \overline{\C}_{+}, a\in P(W)}W_{[n]}^{[a+\Z]}$. It is clear that for $a\in P(V)$ and $b\in P(W)$,
either $a+b\in P(W)$ or $a+b-1\in P(W)$.

Let $(W, Y^{g})$ be a $g$-twisted $V$-module. Then for $v\in V^{[\alpha]}$, we have
$$Y^{g}(v, x)=\sum_{k=0}^{N}\sum_{n \in \alpha + \Z}Y^{g}_{n, k}(v)x^{-n-1}(\log x)^{k}.$$
Let
$$Y^{g}_{k}(v, x)=\sum_{n \in \alpha + \Z}Y^{g}_{n, k}(v)x^{-n-1}.$$
Then
$$Y^{g}(v, x)=\sum_{k=0}^{N}Y^{g}_{k}(v, x)(\log x)^{k}.$$
From equivariance property, it was shown in \cite{B} (cf. Lemma 2.3 \cite{HY}) that
\begin{equation}
Y^g(v, x) = Y^g_0(x^{-\mathcal{N}}v, x),
\end{equation}
where $x^{-\mathcal{N}} = \exp^{-\mathcal{N}\log x}$.

In \cite{B}, Bojko Bakalov derived the following Jacobi type identity for the operator $Y_{0}^{g}$ from an associator formula, which is proved to be equivalent to the duality property in \cite{HY}.
\begin{thm}[\cite{B}]
Let $W$ be a $g$-twisted $V$-module. Then for $u \in V^{[a + \Z]}, v \in V$,  the twisted Jacobi identity
\begin{eqnarray}\label{Jacobi}
x_0^{-1}\delta\left(\frac{x_1 - x_2}{x_0}\right)Y^{g}_{0}(u, x_1)Y^{g}_{0}(v, x_2)
- x_0^{-1}\delta\left(\frac{- x_2 + x_1}{x_0}\right)Y^{g}_{0}(v, x_2)Y^{g}_{0}(u, x_1)\nn
= x_1^{-1}\delta\left(\frac{x_2+x_0}{x_1}\right)\left(\frac{x_2+x_0}{x_1}\right)^{a}
Y^{g}_{0}\left(Y\left(\left(1 + \frac{x_0}{x_2}\right)^{\mathcal{N}}u, x_0\right)v, x_2\right)
\end{eqnarray}
holds.
\end{thm}

As a consequence of Jacobi identity, there are commutator formulas:
\begin{cor}[\cite{B}]
For $u \in V^{[a+\Z]}$, $v\in V^{[b+\Z]}$, $m\in a+\Z$ and $n\in b+\Z$, we have
\begin{eqnarray}\label{commu}
&& {[Y^{g}_{0}(u, x_1), Y^{g}_{0}(v, x_2)] }\nn
&=& \res_{x_0}x_1^{-1}\delta\left(\frac{x_2+x_0}{x_1}\right)
\left(\frac{x_2+x_0}{x_1}\right)^{a} Y^{g}_{0}\left(Y\left(\left(1+\frac{x_0}{x_2}\right)
^{\mathcal{N}}u, x_0\right)v, x_2\right).
\end{eqnarray}
\end{cor}

As the twisted modules for a finite-order automorphism, we also have the weak associativity for the twisted modules for a general automorphism:
\begin{cor}[\cite{HY}]
For $u \in V^{[a+\Z]}$, $v\in V$ and $w \in W$, let $l \in a + \Z$
such that $u(n)w = 0$ for $n \geq l$. Then
\begin{equation}\label{asso}
(x_0+x_2)^l Y^{g}_{0}(u, x_0+x_2)Y^{g}_{0}(v, x_2)w
= (x_2+x_0)^l Y^{g}_{0}\left(Y\left(\left(1+\frac{x_0}{x_2}\right)^{\mathcal{N}}u,
x_0\right)v, x_2\right)w.
\end{equation}
\end{cor}

Finally we have the following equivalence of the main properties of logarithmic $g$-twisted modules:

\begin{thm}[\cite{HY}]\label{jacobi-equiv}
The following properties for a lower bounded logarithmic $g$-twisted
$V$-module are equivalent:
\begin{enumerate}

\item[(1)] The duality property for $Y^{g}$ in Definition \ref{d}.

\item[(2)] The duality property for $Y^{g}_{0}$.

\item[(3)] The property that $Y^{g}_{0}$ is lower truncated
and  the twisted Jacobi identity (\ref{Jacobi}).

\item[(4)] The property that $Y^{g}_{0}$ is lower truncated, the commutator formula (\ref{commu}) and the weak associativity (\ref{asso}).
\end{enumerate}
\end{thm}

\setcounter{equation}{0}
\section{Zhu's algebra for general automorphisms}
In this Section, for a general automorphism $g$ of $V$, we recall the $g$-twisted Zhu's algebra $A_g(V)$ along with functors between categories of certain $V$-modules and $A_g(V)$-modules, constructed in \cite{HY}.

Let $g = \sigma \exp(2\pi i \mathcal{N})$, where $\sigma$ is semisimple and $\mathcal{N}$ is a nilpotent derivation. For homogeneous $u \in V^{[\alpha]}$ and $v \in V$ we define
\begin{equation}
u \circ_g v = \res_x x^{-1-\delta_{\alpha}}Y((1+x)^{\wt u + \mathcal{N} - 1 + \delta_{\alpha} + \alpha}u, x)v,
\end{equation}
where $\delta_{\alpha} = 1$ if $\alpha = 0$, and $\delta_{\alpha} = 0$ if $\alpha \neq 0$. Let $O_g(V)$ be the linear span of all $u\circ_g v$. We define the linear space
\[
A_{g}(V) = V/O_{g}(V)
\]
and a product $*_{g}$ on $V$ as follows:
\begin{eqnarray}\label{d1}
u *_{g} v =  \left\{
\begin{array}{lll}\res_x x^{-1}Y((1+x)^{\wt u + \mathcal{N}}u, x)v && \mathrm{if}\;\; u \in V^{[0]}
\\ 0 &&\mathrm{otherwise}.
\end{array} \right.
\end{eqnarray}

In fact, $A_{g}(V)$ is a quotient of $A(V^{[0]})$ from the following Lemma:
\begin{lemma}\label{l1}
If $\alpha \neq 0$. Then $V^{[\alpha]} \subset O_{g}(V)$.
\end{lemma}

The following theorem states that $A_g(V)$ is an associative algebra with the product $*_{g}$:
\begin{thm}[\cite{HY}]
\begin{enumerate}
\item[(1)]The subspace $O_{g}(V)$ is a two-sided ideal with respect to the product $*_{g}$.
\item[(2)]The quotient space $A_{g}(V)$ is an associative algebra with the product $*_{g}$.
\end{enumerate}
\end{thm}

We shall need the following Lie bracket relation on $A_g(V)$:
\begin{lemma}Assume that $u, v \in V^{[0]}$ are homogeneous. Then
\begin{equation}\label{lierelation}
u*_g v - v*_g u \equiv \res_xY((1+x)^{\wt u +\mathcal{N}-1}u, x)v \;\;\; \mbox{mod}\;\; O_g(V^{[0]}).
\end{equation}
\end{lemma}

For a lower bounded $g$-twisted $V$-module $W = \coprod_{n \in \C, \alpha \in P(W)}W_{[n]}^{[\alpha]}$, let
\begin{equation}
\Omega_g(W) = \{w \in W| u(k)w = 0 \;\mbox{for homogeneous}\; u \in V,\; \mathfrak{R}\{\wt u - k - 1\} < 0\}.
\end{equation}
Then $\Omega_g$ is a functor from the category of lower bounded $g$-twisted $V$-modules to the category of lower bounded $A_g(V)$-modules:
\begin{thm}
The subspace $\Omega_g(W)$ is an $A_g(V)$-module with the action induced by the map $\bar{a} \mapsto o(a)$ for homogeneous $a \in V^{[0]}$ and $\bar{a}$ denotes the image of $a$ in $A_g(V)$, $o(a)$ denotes the operator $Y^g_{\wt a - 1, 0}(a)$ on $W$.
\end{thm}

In \cite{HY}, we also construct a functor $S_g$, which is a right inverse of the functor $\Omega_g$, from the category of lower bounded $A_g(V)$-modules to the category of lower bounded $g$-twisted $V$-modules.

\begin{thm}\label{functor-f}
The functor $\Omega_{g}$ from the category of lower bounded $g$-twisted $V$-modules
to the category of lower bounded $A_{g}(V)$-modules given by $W\mapsto
\Omega_{g}(W)$ has a right inverse, that is, there exists a functor $S_{g}$ from the
category of lower bounded $A_{g}(V)$-modules to the category of
lower bounded $g$-twisted $V$-modules such that $\Omega_{g} \circ S_{g}=1$, where
$1$ is the identity functor on the category of $A_{g}(V)$-modules. Moreover, the functor $S_g$ has the following universal property: Let $M$ be a lower bounded $A_g(V)$-module. For any lower bounded $g$-twisted $V$-module $\tilde{M}$ and an $A_g(V)$-module map $f: M \rightarrow \Omega_g(\tilde{M})$, there exists a unique $g$-twisted module map $\tilde{f}: S_g(M) \rightarrow \tilde{M}$ such that $\tilde{f}|M =  f$.
\end{thm}

For an lower bounded $A_g(V)$-module $M$, let $J_g(M)$ be the maximal proper submodule of $S_g(M)$ (with the property that $J_g(M) \cap M = 0$) and $L_g(M)$ be the irreducible quotient of $S_g(M)$. Then $L_g(M)$ is a lower bounded irreducible $g$-twisted $V$-module satisfying $\Omega_g(L_g(M)) = M$. Thus we have:

\begin{thm}\label{11correspondence}
The functors $L_g$ and $\Omega_g$ give an equivalence between the category of completely reducible lower bounded $g$-twisted $V$-modules and the category of completely reducible lower bounded $A_g(V)$-modules. In particular, $L_g$ and $\Omega_g$ induce a bijection on the isomorphic classes of irreducible lower bounded $g$-twisted $V$-modules and irreducible lower bounded $A_g(V)$-modules.
\end{thm}

\setcounter{equation}{0}
\section{Vertex operator algebras associated to affine Kac-Moody Lie algebras}
In this Section we shall recall the construction of vertex operator algebras associated to affine Kac-Moody Lie algebras, see for instance \cite{LL}.

Let $\mathfrak{g}$ be a finite-dimensional simple Lie algebra over $\C$ with a fixed Cartan subalgebra of $\mathfrak{h}$, and $\Pi = \{\alpha_i\}_{i \in I}$ the set of positive simple roots. We use $\langle\cdot,\cdot\rangle $ to denote the symmetric nondegenerate invariant bilinear form  on $\mathfrak{g}$ such that long roots have square length $2$.

The affine Lie algebra $\widehat{\mathfrak{g}}$ associated with $\mathfrak{g}$ is defined as
\begin{equation*}
\widehat{\mathfrak{g}}=\mathfrak{g}\otimes \C[t,t^{-1}]\oplus \mathbb{C}
\mathbf{k}
\end{equation*}
with $\mathbf{k}$ central and all other bracket relations defined by
\begin{equation*}
 [g\otimes t^m, h\otimes t^n]=[g,h]\otimes t^{m+n}+m\langle
g,h\rangle\delta_{m+n,0}\mathbf{k}
\end{equation*}
for $g,h\in\mathfrak{g}$ and $m,n \in \Z$. We have the decomposition
\begin{equation*}
\widehat{\mathfrak{g}}=\widehat{\mathfrak{g}}_{-} \oplus\widehat{\mathfrak{g}}
_0\oplus\widehat{\mathfrak{g}}_{+}
\end{equation*}
where $\widehat{\mathfrak{g}}_{\pm}=\coprod_{n \in \pm \mathbb{Z}_+} \mathfrak{g}\otimes t^n$ and $\widehat{\mathfrak{g}}_0=\mathfrak{g}\oplus\mathbb{C}\mathbf{k}$.

For a dominant integral weight $\lambda$ of $\mathfrak{g}$ and $\ell\in\mathbb{C}$, we take $M(\lambda,\ell)$ to be the $\widehat{\mathfrak{g}}_0$-module which is the
irreducible $\mathfrak{g}$-module with highest weight $\lambda$ on which
$\mathbf{k}$ acts as $\ell$. We then have the generalized Verma module
\begin{equation*}
V_{\mathfrak{g}}(\lambda, \ell)=U(\widehat{\mathfrak{g}})\otimes_{U(\widehat{\mathfrak{g}}_0\oplus\widehat{\mathfrak{g}}_{+})} M(\lambda,\ell) \cong
U(\widehat{\mathfrak{g}}_-)\otimes_{\C} M(\lambda,\ell),
\end{equation*}
where the linear isomorphism follows from the Poincar${\rm
\acute{e}}$-Birkhoff-Witt theorem.
 The scalar $\ell$ is called the \textit{level} of $V_{\mathfrak{g}}(\lambda, \ell)$. For $g\in\mathfrak{g}$ and $n\in \Z$, we use the notation $g(n)$ to denote the action of $g\otimes t^n$ on a $\widehat{\mathfrak{g}}$-module. Then the
generalized Verma module $V_{\mathfrak{g}}(\lambda, \ell)$ is linearly spanned by vectors of the form
\begin{equation*}
 g_1(-n_1)\cdots g_k(-n_k) u
\end{equation*}
where $g_i \in \mathfrak{g}$, $n_i>0$, and $u\in M(\lambda,\ell)$. Let $J_{\mathfrak{g}}(\lambda, \ell)$ be the maximal submodule of $V_{\mathfrak{g}}(\lambda, \ell)$ and $L_{\mathfrak{g}}(\lambda, \ell) = V_{\mathfrak{g}}(\lambda, \ell)/ J(\lambda, \ell)$. That is, $L(\lambda, \ell)$ is the unique irreducible quotient of $V_{\mathfrak{g}}(\lambda, \ell)$.

We fix a level $\ell\in\C$. When $\ell\neq -h^\vee$, where $h^\vee$ is the dual Coxeter number of $\mathfrak{g}$, the generalized Verma module $V_{\mathfrak{g}}(0, \ell)$ (induced from the one-dimensional $\mathfrak{g}$-module $\mathbb{C}\mathbf{1}$), has the structure of a vertex operator algebra with vacuum $\mathbf{1}$ and vertex operator map determined by
\begin{equation}\label{vertexoperator}
 Y(g(-1)\mathbf{1},x)=g(x)=\sum_{n\in\mathbb{Z}} g(n) x^{-n-1}
\end{equation}
for $g\in\mathfrak{g}$. The conformal vector $\omega$ of $V_{\mathfrak{g}}(0, \ell)$ is
given by
\begin{equation*}
 \omega=\dfrac{1}{2(\ell+h^\vee)}\sum_{i=1}^{\mathrm{dim}\,\mathfrak{g}} u_i(-1)
u_i(-1)\mathbf{1},
\end{equation*}
where $\lbrace u_i\rbrace$ is an orthonormal basis of $\mathfrak{g}$ with
respect to the form $\langle\cdot,\cdot\rangle$.

For any dominant integral weight $\lambda$ of $\mathfrak{g}$, the generalized Verma module
$V_{\mathfrak{g}}(\lambda, \ell)$ is a $V_{\mathfrak{g}}(0, \ell)$-module with vertex operator map defined by (\ref{vertexoperator}). The conformal weight grading on
$V_{\mathfrak{g}}(\lambda, \ell)$ is given by
\begin{equation*}
 \mathrm{wt}\,g_1(-n_1)\cdots g_k(-n_k)u=n_1+\ldots
+n_k+\frac{\langle\lambda,\lambda+2\rho_{\mathfrak{g}}\rangle}{2(\ell+h^\vee)}
\end{equation*}
for $g_i\in\mathfrak{g}$, $n>0$, and $u\in M(\lambda,\ell)$, where
$\rho_{\mathfrak{g}}$ is half the sum of the positive roots of $\mathfrak{g}$.
We observe that $g(n)$ decreases weight by $n$, and so
\begin{equation}\label{l0withgn}
 [L(0), g(n)]=-n g(n)
\end{equation}
for any $g\in\mathfrak{g}$ and $n\in\mathbb{Z}$. In particular, $g(0)$ preserves
weights, so each weight space of $V_{\mathfrak{g}}(\lambda, \ell)$ is a finite-dimensional
$\mathfrak{g}$-module.

The irreducible quotient $L_{\mathfrak{g}}(0, \ell)$ also has a vertex operator algebra structure when $\ell \neq 0$ and $\ell \neq -h^{\vee}$. When $\ell$ is a positive integer, $L_{\mathfrak{g}}(0, \ell)$ is an integrable $\widetilde{\mathfrak{g}}$-module, and the vertex operator algebra $L_{\mathfrak{g}}(0, \ell)$ and its representations give rise to the Wess-Zumino-Novikov-Witten model in the physics literature. It is proved in \cite{FZ} that if $\ell$ is a positive integer, then the vertex operator algebra $L_{\mathfrak{g}}(0, \ell)$ is rational. The set $L_{\mathfrak{g}}(\lambda, \ell)$ with $\lambda \in \mathfrak{h}^{*}$ integral satisfying $\langle \lambda, \theta\rangle \leq \ell$, is a complete list of irreducible $L_{\mathfrak{g}}(0, \ell)$-modules.

\setcounter{equation}{0}
\section{Twisted modules for $V_{\mathfrak{g}}(0, \ell)$}
Let $\mathfrak{g}$ be a finite dimensional Lie algebra and $\mu$ be a Dynkyn diagram automorphism of $\mathfrak{g}$ with order $T$. Then the Lie algebra $\mathfrak{g}$ has the following decomposition:
\begin{equation}
\mathfrak{g} = \coprod_{j =0}^{T-1}\mathfrak{g}^{[j]},
\end{equation}
where
\[
\mathfrak{g}^{[j]} = \{X \in \mathfrak{g}| \mu\cdot X = e^{\frac{2\pi ij}{T}}X\}.
\]
We use $V$ to denote the vertex operator algebras $V_\mathfrak{g}(0, \ell)$ or $L_\mathfrak{g}(0, \ell)$. Then $\mu$ induces an automorphism of $V$, which we still denote by $\mu$. The vertex operator algebra $V$ has the following decomposition:
\begin{equation}
V = \coprod_{j=0}^{T-1}V^{[j]},
\end{equation}
where
\[
V^{[j]} = \{v \in V| \mu \cdot v = e^{\frac{2\pi ij}{T}}v \}.
\]

Let $\varphi$ be an arbitrary element in  $Aut(\mathfrak{g})$. Then $\varphi = \mu\exp(\mbox{ad}\; h)\exp(\mbox{ad}\; e)$, where $e$ is a nilpotent element in $\mathfrak{g}^{[0]}$ and $h$ is a semisimple element in $\mathfrak{g}$. In the remaining of this paper, we will let $\sigma = \mu\exp(\mbox{ad}\; e)$. It was shown in \cite{L} that there is a one-to-one correspondence between $\varphi$-twisted $V$-modules and $\sigma$-twisted $V$-modules. In this Section, by determining the $\sigma$-twisted Zhu's algebra $A_{\sigma}(V)$, we classify the irreducible lower bounded $\sigma$-twisted $V$-modules and determine the complete reducibility for these modules. We shall need the following definition:
\begin{defn}{\rm
For $k \in \Z_{+}$ and $v \in V_{\mathfrak{g}}(0, \ell)$, we say that $v$ has $({\rm ad}\; e)$-block size $k$ if $({\rm ad}\; e)^{k-1}v \neq 0$ and $({\rm ad}\; e)^{k}v = 0$.}
\end{defn}

We first show that
\begin{lemma}\label{genlemma1}
The elements in $A_{\sigma}(V)$ are of the form
\begin{equation}\label{standardform1}
[g_1(-n_1)\cdots g_k(-n_k){\bf 1}],
\end{equation}
where $g_i \in \mathfrak{g}^{[0]}$, $n_i \in \Z_{+}$ and $k \in \N$.
\end{lemma}
\pf By Lemma $\ref{l1}$, we know that $V^{[j]} \subset O_{\mu}(V)$ for $j \neq 0$, it remains to prove the lemma for $v \in V^{[0]}$. Let $v = h_1(-i_1)\cdots h_n(-i_n){\bf 1}$ for $h_i \in \mathfrak{g}^{[\alpha_i]}$ and $\sum_{i = 1}^n\alpha_i = 0$ mod $T$. We prove by induction on $n$ that $[v]$ is of the form (\ref{standardform1}). If $\alpha_i$'s are $0$ for all $i = 1, \dots, n$, then we are done. Otherwise, we can re-arrange the order of $h_i$ such that there exists $1 < m \leq n$ such that $\alpha_i \neq 0$ for $i = 1, \dots, m$ and $\alpha_i = 0$ for $i = m+1, \dots, n$. Since for arbitrary $k \in \N$,
\begin{eqnarray*}
&&\res_x x^{-1-k}(1+x)^{\alpha_1 + \mbox{ad}\; e}Y(h_1(-1){\bf 1}, x)v\nn
& = & \sum_{i \geq 0}\bigg(\binom{\alpha_1+ \mbox{ad}\; e}{i}h_1\bigg)(i-k-1)v \in O_{\sigma}(V),
\end{eqnarray*}
we have that
\begin{eqnarray*}
h_1(-k-1)v \equiv -\sum_{i \geq 1}\bigg(\binom{\alpha_1+ \mbox{ad}\; e}{i}h_1\bigg)(i-k-1)v \;\;\; {\rm mod}\; O_{\sigma}(V)
\end{eqnarray*}
and thus
\begin{eqnarray*}
&& h_1(-i_1)\cdots h_n(-i_n){\bf 1}\nn
&\equiv & \sum h_2'(-i_2')\cdots h_n'(-i_n'){\bf 1} \;\;\; {\rm mod}\; O_{\sigma}(V),
\end{eqnarray*}
where the sum ranges over $h_i' \in \mathfrak{g}^{[\alpha_i']}$ and $\sum_{i = 2}^n\alpha_i' = 0$ mod $T$. Then the lemma follows from induction. \epfv

The associative algebra $A_{\sigma}(V_{\mathfrak{g}}(0, \ell))$ is given by the following theorem.
\begin{thm}\label{unithm}
The associative algebra $A_{\sigma}(V_{\mathfrak{g}}(0, \ell))$ is canonically isomorphic to $U(\mathfrak{g}^{[0]})$.
\end{thm}
\pf
Let $[b] \in A_{\sigma}(V_{\mathfrak{g}}(0, \ell))$ denote the image of $b \in V_{\mathfrak{g}}(0, \ell)$.
We define a linear map $i: \mathfrak{g}^{[0]} \longrightarrow A_{\sigma}(V_{\mathfrak{g}}(0, \ell))$ by
\[
i: g \longmapsto [g(-1){\bf 1} + \ell \langle e, g\rangle{\bf 1}]
\]
for $g \in \mathfrak{g}^{[0]}$. From (\ref{lierelation}), we have for $a, b \in \mathfrak{g}^{[0]}$,
\begin{eqnarray}\label{com}
&&[a(-1){\bf 1}, b(-1){\bf 1}] \nn
&=& (a(-1){\bf 1})*_g (b(-1){\bf 1}) - (b(-1){\bf 1})*_g (a(-1){\bf 1})\nn
&=& [a, b](-1){\bf 1} + \ell \langle ({\rm ad}\;e)a, b\rangle{\bf 1}\nn
&=& [a, b](-1){\bf 1} + \ell \langle [e, a], b\rangle{\bf 1}.
\end{eqnarray}
Using (\ref{com}),
\begin{eqnarray*}
&&[i(a), i(b)]\nn
&=& [[a(-1){\bf 1} + \ell \langle e, a\rangle{\bf 1}, b(-1){\bf 1} + \ell \langle e, b\rangle{\bf 1}]]\nn
&=& [[a(-1){\bf 1}, b(-1){\bf 1}]]\nn
&=& [[a, b](-1){\bf 1} + \ell \langle [e, a], b\rangle{\bf 1}]\nn
&=& [[a, b](-1){\bf 1} + \ell \langle e, [a, b]\rangle{\bf 1}]\nn
&=& i([a, b]).
\end{eqnarray*}
Thus $i$ induces an algebra homomorphism $I: U(\mathfrak{g}^{[0]})\longrightarrow A_{\sigma}(V_{\mathfrak{g}}(0, \ell))$ by
\begin{equation}\label{mapi}
g_1\cdots g_n \longmapsto i(g_1)*_{\sigma} \cdots *_{\sigma} i(g_n).
\end{equation}
By the definition of $*_{\sigma}$, we have
\[
[(g_1(-1){\bf 1})*_{\sigma}\cdots *_{\sigma}(g_n(-1){\bf 1})] \in [g_1(-1)\cdots g_n(-1){\bf 1}] + R_n,
\]
where $R_n$ is the linear span of elements of the form $[\tilde{g}_1(-1)\cdots \tilde{g}_k(-1){\bf 1}]$ for $\tilde{g}_i \in \mathfrak{g}^{[0]}$ and $0 \leq k < n$. Thus by induction on $n$, the elements $[g_1(-1)\cdots g_n(-1){\bf 1}]$ are in $I(U(\mathfrak{g}^{[0]}))$.

From Lemma \ref{genlemma1}, we know that the elements in $A_{\sigma}(V_{\mathfrak{g}}(0, \ell))$  are of the form
\begin{equation*}
[g_1(-i_1-1)\cdots g_k(-i_k-1){\bf 1}],
\end{equation*}
where $g_m \in \mathfrak{g}^{[0]}$, $i_m, k \in \N$.
Since for $a \in \mathfrak{g}^{[0]}, b \in V^{[0]}$ and $n \geq 0$,
\begin{eqnarray*}
&&\res_x Y\bigg(\frac{(1+x)^{1 + {\rm ad}\; e}}{x^{2+n}}(a(-1){\bf 1}), x\bigg)b\nn
&=& \sum_{i \geq 0}\bigg(\binom{1+ {\rm ad}\; e}{i}a\bigg)(i-n-2)b\nn
&=& a(-n-2)b + ((1 + {\rm ad}\; e)a)(-n-1)b + \cdots\nn
& \in & O_{\sigma}(V_{\mathfrak{g}}(\ell, 0)),
\end{eqnarray*}
we have
\[
[g_1(-i_1-1)\cdots g_k(-i_k-1){\bf 1}] \in R_{k+1},
\]
which is in the image of $I$; thus $I$ is surjective.

We still need to prove the injectivity. In fact, the elements of $O_{\sigma}(V_{\mathfrak{g}}(0, \ell))$ are spanned by
\begin{eqnarray*}
u \circ_{\sigma} v &= & Y\bigg(\frac{(1+x)^{\wt u +\alpha + \delta_{\alpha} - 1+ {\rm ad}\; e}}{x^{1+\delta_{\alpha}}}u, x\bigg)v\nn
&= & Y\bigg(\frac{(1+x)^{\wt u +\alpha + \delta_{\alpha} - 1}}{x^{1+ \delta_{\alpha}}}u, x\bigg)v + Y\bigg(\frac{(1+x)^{\wt u+\alpha + \delta_{\alpha} - 1}((1+x)^{{\rm ad}\; e}-1)}{x^2}u, x\bigg)v\nn
&= & u\circ_{\mu} v + Y\bigg(\frac{(1+x)^{\wt u+\alpha + \delta_{\alpha} - 1}((1+x)^{{\rm ad}\; e}-1)}{x^2}u, x\bigg)v.
\end{eqnarray*}
It is obvious that each summand of $u \circ_{\mu} v$ and each summand of
\[
Y\bigg(\frac{(1+x)^{\wt u+\alpha + \delta_{\alpha} - 1}((1+x)^{{\rm ad}\; e}-1)}{x^2}u, x\bigg)v
\]
have different $({\rm ad}\; e)$-block sizes. Thus they are linearly independent. From the proof of the surjectivity, $I(U(\mathfrak{g}^{[0]}))$ is spanned by elements of the form
\[
[g_1(-1)\cdots g_n(-1){\bf 1}],  \; g_i \in \mathfrak{g}^{[0]},
\]
it suffices to show $O_{\mu}(V_{\mathfrak{g}}(0, \ell))$ is spanned by elements of the form
\[
(g(-n-2)+g(-n-1))v,\; g \in \mathfrak{g}^{[0]}, v \in V_{\mathfrak{g}}(0, \ell)
\]
and elements of the form
\[
g_1(-i_1-1)\cdots g_k(-i_k-1)v,\; g_j \in \mathfrak{g}^{[\alpha_i]}, v \in V_{\mathfrak{g}}(0, \ell)
\]
for some $\alpha_i \neq 0$ and $i_j, k \in \N$.

Let $u = g_1(-i_1-1)\cdots g_n(-i_n-1){\bf 1} \in V_{\mathfrak{g}}(0, \ell)$ for $h_j \in \mathfrak{g}^{[\alpha_i]}$ and $i_j \in \N$, we define the length of $u$ to be $n$.

If $\alpha_i's$ are all $0$, then by the result of the untwisted Zhu's algebra (see \cite{FZ}), $u\circ_{\mu} v$ is spanned by elements of the form
\[
(g(-n-2)+g(-n-1))v, g \in \mathfrak{g}^{[0]}, v \in V_{\mathfrak{g}}(0, \ell).
\]
Otherwise, we can re-arrange the order of $h_i$'s such that there exists $1 < m \leq n$ such that $\alpha_i \neq 0$ for $i = 1, \dots, m$ and $\alpha_i = 0$ for $i = m+1, \dots, n$. Since the highest weight summand of $u \circ_{\mu} v$ is
\[
\res_x x^{-1-\delta_{\alpha}}Y(u, x)v
\]
for $u \in V^{[\alpha]}$ and
\[
\res_x x^{-2}Y(u, x)v = \res_{x} x^{-1}Y(L(-1)u, x)v,
\]
without loss of generality, we assume the highest weight summand of $u \circ_{\mu} v$ is of the form $\res_x  x^{-1}Y(u, x)v$, which is
\[
 u_{-1}v = g_1(-i_1-1)\cdots g_n(-i_n-1)v + \;\mbox{shorter length terms}.
\]
Thus $I$ is injective.
\epfv

In \cite{B}, Bojko Bakalov defined a {\it $\sigma$-twisted affinization $\tilde{\mathfrak{g}}[\sigma]$} of $\mathfrak{g}$, which consists of elements of the form $a \otimes t^m$ for $a \in \mathfrak{g}^{[i]}$, $m \in \frac{i}{T} + \Z$ with Lie bracket defined by
\[
[a \otimes t^m, b \otimes t^n] = [a, b]\otimes t^{m+n} + \delta_{m+n, 0}(m\langle a, b \rangle + \langle e, [a, b]\rangle){\bf k}
\]
for $a \in \mathfrak{g}^{[i]}, m \in \frac{i}{T} + \Z, b \in \mathfrak{g}^{[j]}, n \in \frac{j}{T} + \Z$. He showed that the category of lower bounded $\sigma$-twisted $V_{\mathfrak{g}}(0, \ell)$-modules are equivalent to the category of restricted $\tilde{\mathfrak{g}}[\sigma]$-modules. In this Section, we will determine these modules explicitly using Theorem \ref{11correspondence} and Theorem \ref{unithm}.

Let $\tilde{\mathfrak{g}}[\mu]$ denote the twisted affine Lie algebra associated with the Lie algebra $\mathfrak{g}$ and the diagram automorphism $\mu$. We write the elements of $\tilde{\mathfrak{g}}[\mu]$ in the form $a(m)$ for $a \in \mathfrak{g}^{[i]}, m \in \frac{i}{T} + \Z$. Then $\tilde{\mathfrak{g}}[\mu]$ is isomorphic to $\tilde{\mathfrak{g}}[\sigma]$ as Lie algebras under the canonical map:
\begin{eqnarray}\label{affineiso}
\phi: \tilde{\mathfrak{g}}[\sigma] &\longrightarrow & \tilde{\mathfrak{g}}[\mu]\nn
a \otimes t^m &\longmapsto & a(m) - \delta_{m, 0}\langle e, a\rangle{\bf k}.
\end{eqnarray}
Recall that $\tilde{\mathfrak{g}}[\mu]$ has the canonical triangular decomposition
\[
\tilde{\mathfrak{g}}[\mu] = \tilde{\mathfrak{g}}[\mu]_{+} \oplus \tilde{\mathfrak{g}}[\mu]_0 \oplus \tilde{\mathfrak{g}}[\mu]_{-},
\]
where $\tilde{\mathfrak{g}}[\mu]_{+}$ (resp. $\tilde{\mathfrak{g}}[\mu]_{-}$) is the Lie subalgebra of $\tilde{\mathfrak{g}}[\mu]$ consisting of elements with positive (resp. negative) weights and $\tilde{\mathfrak{g}}[\mu]_{0} = \mathfrak{g}^{[0]} \oplus \C{\bf k}$.
Let $\mathfrak{h}^{[0]}$ be a Cartan subalgebra of $\mathfrak{g}^{[0]}$. Fix a dominant integral weight $\lambda \in (\mathfrak{h}^{[0]})^*$, let $M(\lambda, \ell)$ be a  $\tilde{\mathfrak{g}}[\mu]_{0}$-module which is  the irreducible $\mathfrak{g}^{[0]}$-module with highest weight $\lambda$ on which ${\bf k}$ acts as $\ell$. We then have the generalized Verma module associated with $M(\lambda, \ell)$:
\[
V_{(\mathfrak{g}, \mu)}(\lambda, \ell): = U(\tilde{\mathfrak{g}}[\mu]) \otimes_{\tilde{\mathfrak{g}}[\mu]_{+} \oplus \tilde{\mathfrak{g}}[\mu]_0} M(\lambda, \ell).
\]
Let $J_{(\mathfrak{g}, \mu)}(\lambda, \ell)$ be the maximal proper $\tilde{\mathfrak{g}}[\mu]$-submodule of $V_{(\mathfrak{g}, \mu)}(\lambda, \ell)$. Let
\[
L_{(\mathfrak{g}, \mu)}(\lambda, \ell) =V_{(\mathfrak{g}, \mu)}(\lambda, \ell)/J_{(\mathfrak{g}, \mu)}(\lambda, \ell)
\]
be the irreducible quotient of $V_{(\mathfrak{g}, \mu)}(\lambda, \ell)$. Then $V_{(\mathfrak{g}, \mu)}(\lambda, \ell)$ and $L_{(\mathfrak{g}, \mu)}(\lambda, \ell)$ can be viewed as $\tilde{\mathfrak{g}}[\sigma]$-modules through the Lie algebra isomorphism (\ref{affineiso}), and hence they are also $\sigma$-twisted $V_{\mathfrak{g}}(0, \ell)$-modules.

\begin{prop}
The generalized Verma module $V_{(\mathfrak{g}, \mu)}(\lambda, \ell)$ is the same as the $\sigma$-twisted module $S_{\sigma}(M(\lambda, \ell))$ constructed in \cite{HY} (see Section $3$), and every lower bounded $\sigma$-twisted $V_{\mathfrak{g}}(0, \ell)$-module $W$ such that $\Omega_{\sigma}(W) = M(\lambda, \ell)$ is a quotient of $V_{(\mathfrak{g}, \mu)}(\lambda, \ell)$.
\end{prop}
\pf Since $A_{\sigma}(V_{\mathfrak{g}(0, \ell)}) = U(\mathfrak{g}^{[0]})$, the $\sigma$-twisted modules $V_{(\mathfrak{g}, \mu)}(\lambda, \ell)$ also satisfies the universal property for $S_{\sigma}(M(\lambda, \ell))$ in Theorem \ref{functor-f}. \epfv

As a consequence, we provide the list of all lower bounded irreducible $\sigma$-twisted $V_{\mathfrak{g}}(0, \ell)$-modules:
\begin{cor}
Let $\ell$ be a positive integer and $\lambda$ be a dominant integral weight for $\mathfrak{g}^{[0]}$, then $L_{(\mathfrak{g}, \mu)}(\lambda, \ell)$ is a lower bounded irreducible $\sigma$-twisted module for the vertex operator algebra $V_{\mathfrak{g}}(0, \ell)$. Conversely, every lower bounded irreducible $\sigma$-twisted $V_{\mathfrak{g}}(0, \ell)$-module of level $\ell$ is of the form  $L_{(\mathfrak{g}, \mu)}(\lambda, \ell)$ for some dominant integral weight $\lambda$.
\end{cor}

\setcounter{equation}{0}
\section{Twisted modules for $L_{\mathfrak{g}}(0, \ell)$}
In the Section, we will prove the complete reducibility for the lower bounded $\sigma$-twisted $L_{\mathfrak{g}}(0, \ell)$-modules and describe the irreducible $\sigma$-twisted modules, using the $\sigma$-twisted Zhu's algebra for $L_{\mathfrak{g}}(0, \ell)$.

Fix a Cartan subalgebra $\mathfrak{h}$ of $\mathfrak{g}$ and a base for the root system, let $\theta$ be the highest root. Fix nonzero vectors $e_{\theta} \in \mathfrak{g}_{\theta}$ and $f_{\theta} \in \mathfrak{g}_{-\theta}$ such that $\langle e_{\theta}, f_{\theta}\rangle = 1$.

The maximal submodule $J_{\mathfrak{g}}(0, \ell)$ of $V_{\mathfrak{g}}(0, \ell)$ is generated by the element $e_{\theta}(-1)^{\ell + 1}{\bf 1}$. Let $\mu$ be the nontrivial diagram automorphism of $\mathfrak{g}$. Then in the case $\mathfrak{g} = A_{2n-1}, D_{n}$ or $E_6$ for $n \geq 3$, we know $e_{\theta} \in \mathfrak{g}^{[0]}$ and in the case $\mathfrak{g} = A_{2n}$ for $n \geq 1$, $e_{\theta} \in \mathfrak{g}^{[1]}$ (see \cite{K}). We first determine the image of  $e_{\theta}(-1)^{\ell+1}{\bf 1}$ in $A_{\sigma}(V_{\mathfrak{g}}(0, \ell))$ in the following lemma:
\begin{lemma}\label{mainlemma}
In $A_{\sigma}(V_{\mathfrak{g}}(0, \ell))$, we have for $k \in \N$,
\begin{equation}
([e_{\theta}(-1){\bf 1} + \ell\langle e, e_{\theta}\rangle{\bf 1}])^{k} = [e_{\theta}(-1)^k{\bf 1}] + \sum_{i=0}^{k-1}\binom{k}{i}\binom{\ell - i}{k-i}(k-i)![e_{\theta}(-1)^i{\bf 1}].
\end{equation}
In particular,
\begin{equation}
([e_{\theta}(-1){\bf 1} + \ell\langle e, e_{\theta}\rangle{\bf 1}])^{\ell+1} = [e_{\theta}(-1)^{\ell+1}{\bf 1}].
\end{equation}
\end{lemma}
\pf
If $e_{\theta} \in \mathfrak{g}^{[1]}$, then the lemma obviously holds since both sides are $0$. We assume $e_{\theta} \in \mathfrak{g}^{[0]}$,
it suffices to prove the case when the nilpotent element $e = f_{\theta}$. We can show by induction on $k$ that
\begin{equation}\label{little1}
f_{\theta}(1)e_{\theta}(-1)^k{\bf 1} = k(\ell-k+1)e_{\theta}(-1)^{k-1}{\bf 1}.
\end{equation}
Also it is easy to see the following formulas
\begin{equation}\label{little2}
f_{\theta}(i)e_{\theta}(-1)^k{\bf 1} = 0\;\;\; \mbox{for}\; i \geq 2,
\end{equation}
\begin{equation}\label{little3}
h_{\theta}(i)e_{\theta}(-1)^k{\bf 1} = 0\;\;\; \mbox{for}\; i \geq 1.
\end{equation}
Using formulas (\ref{little1}), (\ref{little2}) and (\ref{little3}), we have
\begin{equation}\label{little4}
(e_{\theta}(-1){\bf 1}) *_{\sigma} (e_{\theta}(-1))^k{\bf 1} = (e_{\theta}(-1))^{k+1}{\bf 1} - 2k (e_{\theta}(-1))^k{\bf 1} - k(\ell - k+1)(e_{\theta}(-1))^{k-1}{\bf 1}.
\end{equation}
Now we prove the lemma by induction on $k$. We assume that
\begin{equation*}
([e_{\theta}(-1){\bf 1} + \ell{\bf 1}])^{k} = [e_{\theta}(-1)^k{\bf 1}] + \sum_{i=0}^{k-1}\binom{k}{i}\binom{\ell - i}{k-i}(k-i)![e_{\theta}(-1)^i{\bf 1}].
\end{equation*}
Then
\begin{eqnarray*}
&&([e_{\theta}(-1){\bf 1} + \ell{\bf 1}])^{k+1}\nn
&=&([e_{\theta}(-1){\bf 1} + \ell{\bf 1}])*_{\sigma}([e_{\theta}(-1){\bf 1} + \ell{\bf 1}])^{k}\nn
&=& ([e_{\theta}(-1){\bf 1} + \ell{\bf 1}])*_{\sigma} [e_{\theta}(-1)^k{\bf 1}] + \sum_{i=0}^{k-1}\binom{k}{i}\binom{\ell - i}{k-i}(k-i)!([e_{\theta}(-1){\bf 1} + \ell{\bf 1}])*_{\sigma}[e_{\theta}(-1)^i{\bf 1}]\nn
&=& [e_{\theta}(-1)^{k+1}{\bf 1}] + \bigg(k\binom{\ell - k +1}{1}-2k + \ell\bigg)[e_{\theta}(-1)^k{\bf 1}]\nn
&& +\; \sum_{i=1}^{k-1}\bigg(\binom{k}{i-1}\binom{\ell -i+1}{k-i+1}(k-i+1)! - \binom{k}{i}\binom{\ell - i}{k-i}(k-i)!2i\nn
&& -\; \binom{k}{i+1}\binom{\ell-i-1}{k-i-1}(k-i-1)!(i+1)(\ell-i) + \binom{k}{i}\binom{\ell-i}{k-i}(k-i)!\ell\bigg)[e_{\theta}(-1)^i{\bf 1}]\nn
&& + \; \bigg(-\binom{k}{1}\binom{\ell-1}{k-1}(k-1)!\ell + \binom{\ell}{k}k!\ell\bigg)[{\bf 1}]\nn
&=& [e_{\theta}(-1)^{k+1}{\bf 1}] + \sum_{i=0}^{k}\binom{k+1}{i}\binom{\ell - i}{k-i+1}(k-i+1)![e_{\theta}(-1)^i{\bf 1}],
\end{eqnarray*}
where the last step follows from straightforward computations. \epfv

\begin{rema}{\rm Let $\theta^0$ be the highest root in the Lie subalgebra $\mathfrak{g}^{[0]}$ and let $e_{\theta^0}$ be the root vector corresponding to $\theta^0$. Then the result in Lemma \ref{mainlemma} also holds if we replace $e_{\theta}$ by $e_{\theta^0}$.}
\end{rema}

The following theorem describes the $\sigma$-twisted Zhu's algebra for $L_{\mathfrak{g}}(0, \ell)$:
\begin{thm}\label{maintheorem1}Let $\ell$ be a positive integer and let $\mathfrak{g}$ be a finite dimensional simple Lie algebra and $\mu$ be a diagram automorphism of $\mathfrak{g}$. Let $\sigma = \mu \exp(\mbox{ad}\; e)$ for a nilpotent element $e\in \mathfrak{g}^{[0]}$.
\begin{itemize}
\item[(1)]
If $\mu$ is trivial, i.e., $\sigma$ is an inner automorphism. Then the $\sigma$-twisted Zhu's algebra $A_{\sigma}(L_{\mathfrak{g}}(0, \ell))$ is canonically isomorphic to $U(\mathfrak{g})/\langle e_{\theta}^{\ell + 1}\rangle$.
\item[(2)]
If $\mu$ is nontrivial and $\mathfrak{g}$ is not of type $A_{2n}$ for $n \geq 1$. Then the $\sigma$-twisted Zhu's algebra $A_{\sigma}(L_{\mathfrak{g}}(0, \ell))$ is a quotient of $U(\mathfrak{g}^{[0]})/\langle e_{\theta}^{\ell + 1}\rangle$.
\item[(3)]
If $\mathfrak{g}$ is of type $A_{2n}$ for $n \geq 1$ and $\mu$ has order $2$. Then the $\sigma$-twisted Zhu's algebra $A_{\sigma}(L_{\mathfrak{g}}(0, \ell))$ is a quotient of $U(\mathfrak{g}^{[0]})/\langle e_{\theta^0}^{\ell + 1}\rangle$, where $e_{\theta^0}$ is the root vector corresponding to the highest root $\theta^0$ of the Lie subalgebra $\mathfrak{g}^{[0]}$.
\end{itemize}
\end{thm}
\pf
Similar to Proposition $1.4.2$ in \cite{FZ}, it is easy to see $A_{\sigma}(J_{\mathfrak{g}}(0, \ell))$ is a two-sided ideal in $A_{\sigma}(V_{\mathfrak{g}}(0, \ell))$ and $A_{\sigma}(L_{\mathfrak{g}}(0, \ell)) = A_{\sigma}(V_{\mathfrak{g}}(0, \ell))/A_{\sigma}(J_{\mathfrak{g}}(0, \ell))$. From Lemma \ref{mainlemma}, we know $[e_{\theta}(-1)^{\ell + 1}{\bf 1}] \in A_{\sigma}(V_{\mathfrak{g}}(0, \ell))$ (or $[e_{\theta^0}(-1)^{\ell + 1}{\bf 1}]$) are mapped to $e_{\theta}^{\ell + 1} \in U(\mathfrak{g}^{[0]})$ (or $e_{\theta^0}^{\ell + 1}$) under the canonical isomorphism $I^{-1}$. Thus the two-sided ideal $\langle e_{\theta}^{\ell + 1} \rangle$ (or $\langle e_{\theta^0}^{\ell + 1} \rangle$) $\subset A_{\sigma}(J_{\mathfrak{g}}(0, \ell))$, this proves $(2)$ and $(3)$.

We still need to prove $(1)$. Note that every element in $J_{\mathfrak{g}}(0, \ell)$ can be written as a linear combination of elements of type
\begin{equation}\label{elementsm}
a_1(-i_1)\cdots a_n(-i_n)e_{\theta}(-1)^{\ell + 1}{\bf 1}
\end{equation}
for $a_m \in \mathfrak{g}$ and $i_m \geq 0$, we prove by induction on $n$ that
\[
[a_1(-i_1)\cdots a_n(-i_n)e_{\theta}(-1)^{\ell + 1}{\bf 1}] \in \langle [e_{\theta}(-1)^{\ell + 1}{\bf 1}]\rangle.
\]
Since
\begin{eqnarray*}
&&\res_x Y\bigg(\frac{(1+x)^{1 + {\rm ad}\; e}}{x^{2+n}}(a(-1){\bf 1}), x\bigg)b\nn
&=& \sum_{i \geq 0}\bigg(\binom{1+ {\rm ad}\; e}{i}a\bigg)(i-n-2)b\nn
&\in &  O_{\sigma}(V_{\mathfrak{g}}(0, \ell))
\end{eqnarray*}
for $a \in \mathfrak{g}$, $b \in V$ and $n \geq 0$, and also by induction hypothesis, it suffices to show
\begin{equation}\label{indhyp1}
[a_1(-1)a_2(-i_2)\cdots a_n(-i_n)e_{\theta}(-1)^{\ell + 1}{\bf 1}] \in \langle [e_{\theta}(-1)^{\ell + 1}{\bf 1}]\rangle
\end{equation}
and
\begin{equation}\label{indhyp2}
[a_1(0)a_2(-i_2)\cdots a_n(-i_n)e_{\theta}(-1)^{\ell + 1}{\bf 1}] \in \langle [e_{\theta}(-1)^{\ell + 1}{\bf 1}]\rangle.
\end{equation}
In fact,
\begin{eqnarray*}
&&(a_1(-1){\bf 1})*_{\sigma} (a_2(-i_2)\cdots a_n(-i_n)e_{\theta}(-1)^{\ell + 1}{\bf 1}) - (a_2(-i_2)\cdots a_n(-i_n)e_{\theta}(-1)^{\ell + 1}{\bf 1}) *_{\sigma} (a_1(-1){\bf 1})\nn
&=& \sum_{i\in \N}\bigg(\binom{{\rm ad}\; e}{i}a_1\bigg)(i)(a_2(-i_2)\cdots a_n(-i_n)e_{\theta}(-1)^{\ell + 1}{\bf 1})\nn
&=& a_1(0)a_2(-i_2)\cdots a_n(-i_n)e_{\theta}(-1)^{\ell + 1}{\bf 1} + \sum_{i\in \Z_{+}}\bigg(\binom{{\rm ad}\; e}{i}a_1\bigg)(i)(a_2(-i_2)\cdots a_n(-i_n)e_{\theta}(-1)^{\ell + 1}{\bf 1}),
\end{eqnarray*}
note that the left-hand-side and the second term in the right-hand-side are in $\langle  [e_{\theta}(-1)^{\ell + 1}{\bf 1}] \rangle$ by induction hypothesis, (\ref{indhyp2}) follows. (\ref{indhyp1}) follows from induction hypothesis, (\ref{indhyp2}) and
\begin{eqnarray*}
&&(a_1(-1){\bf 1})*_{\sigma} (a_2(-i_2)\cdots a_n(-i_n)e_{\theta}(-1)^{\ell + 1}{\bf 1})\nn
&=& \sum_{i\in \N}\bigg(\binom{{\rm ad}\; e + 1}{i}a_1\bigg)(i-1)(a_2(-i_2)\cdots a_n(-i_n)e_{\theta}(-1)^{\ell + 1}{\bf 1})\nn
&=& a_1(-1)a_2(-i_2)\cdots a_n(-i_n)e_{\theta}(-1)^{\ell + 1}{\bf 1} \nn
&& + \sum_{i\in \Z_{+}}\bigg(\binom{{\rm ad}\; e + 1}{i}a_1\bigg)(i-1)(a_2(-i_2)\cdots a_n(-i_n)e_{\theta}(-1)^{\ell + 1}{\bf 1}).
\end{eqnarray*}\epfv

As a consequence, we describe the lower bounded $\sigma$-twisted modules for $L_{\mathfrak{g}}(0, \ell)$:
\begin{thm}\label{mainthmpaper}Let $\ell$ be a positive integer and $\mathfrak{g}$ be a finite dimensional simple Lie algebra.
\begin{itemize}
\item[(1)]If $\mu$ is trivial. i.e., $\sigma$ is an inner automorphism. Then the set of lower bounded irreducible $\sigma$-twisted modules of level $\ell$ for $L_{\mathfrak{g}}(0, \ell)$ consists of $L_{\mathfrak{g}}(\lambda, \ell)$ with dominant integral weight satisfying $\langle \lambda, \theta \rangle \leq \ell$. In particular, there are only finitely many irreducible grading-restricted $\sigma$-twisted modules for $L_{\mathfrak{g}}(0, \ell)$.
\item[(2)]If $\mu$ is nontrivial and $\mathfrak{g}$ is not of type $A_{2n}$ for $n \geq 1$. Then every lower bounded irreducible $\sigma$-twisted modules of level $\ell$ for $L_{\mathfrak{g}}(0, \ell)$ is of the form $L_{(\mathfrak{g}, \mu)}(\lambda, \ell)$ with dominant integral weight satisfying $\langle \lambda, \theta \rangle \leq \ell$.
\item[(3)]If $\mathfrak{g}$ is of type $A_{2n}$ for $n \geq 1$ and $\mu$ has order $2$. Then every lower bounded irreducible $\sigma$-twisted modules of level $\ell$ for the vertex operator algebra $L_{\mathfrak{g}}(0, \ell)$ is of the form $L_{(\mathfrak{g}, \mu)}(\lambda, \ell)$ with dominant integral weight satisfying $\langle \lambda, \theta^0 \rangle \leq \ell$.
\item[(4)]The lower bounded $\sigma$-twisted $L_{\mathfrak{g}}(0, \ell)$-modules are completely reducible.
\end{itemize}
\end{thm}

\pf Let $W$ be a $\sigma$-twisted $L_{\mathfrak{g}}(0, \ell)$-module. Then $\Omega_{\sigma}(W)$ is an $A_{\sigma}(L_{\mathfrak{g}}(0, \ell))$-module. By Theorem \ref{maintheorem1}, $\Omega_{\sigma}(W)$ is a $\mathfrak{g}^{[0]}$-module satisfying that $e_{\theta}^{\ell + 1} = 0$ (or $e_{\theta^0}^{\ell + 1} = 0$). Thus
$\Omega_{\sigma}(W)$ is a direct sum of irreducible $\mathfrak{g}^{[0]}$-modules and hence $W$ is completely reducible. \epfv

We still cannot determine whether all the irreducible $\sigma$-twisted $V_{\mathfrak{g}}(0, \ell)$-modules in the lists in Theorem \ref{mainthmpaper} are irreducible $\sigma$-twisted $L_{\mathfrak{g}}(0, \ell)$-modules when $\mu$ is nontrivial. It depends on a precise characterization of the $\sigma$-twisted Zhu's algebra. However, we will show that the module category of lower bounded $\sigma$-twisted modules is equivalent to the module category of lower bounded $\mu$-twisted modules. We need the following lemma first:
\begin{lemma}\label{exception}Let $f$ be a root vector of $\mathfrak{g}$ contained in $\mathfrak{g}^{[i]}$ for $i \neq 0$. Then
\[
Y^{\sigma}_0(f(-1)^{\ell + 1}{\bf 1}, x) = Y^{\sigma}_0(f(-1){\bf 1}, x)^{\ell + 1}.
\]
\end{lemma}
\pf It is easy to see from the commutator formula (\ref{commu}) that the operators $Y^{\sigma}_0(f(-1){\bf 1}, x)$
and $Y^{\sigma}_0(f(-1)^{n}{\bf 1}, x)$ commute for each $n \in \N$. Then the lemma follows from
\[
Y^{\sigma}_0(f(-1)^{n + 1}{\bf 1}, x) = Y^{\sigma}_0(f(-1){\bf 1}, x)Y^{\sigma}_0(f(-1)^{n}{\bf 1}, x).
\]\epfv

\begin{thm}\label{ind}
The following module categories are equivalent:
\begin{itemize}
\item[(1)]Lower bounded $\sigma$-twisted $L_{\mathfrak{g}}(0, \ell)$-modules;
\item[(2)]Lower bounded $\mu$-twisted $L_{\mathfrak{g}}(0, \ell)$-modules;
\item[(3)]Direct sums of standard $\tilde{\mathfrak{g}}[\mu]$-modules of level $\ell$.
\end{itemize}
\end{thm}
\pf It suffices to show that the simple objects in these module categories have one-to-one bijections. $(2) \Leftrightarrow(3)$ is Proposition 5.6 in \cite{L}. For $(3)\Rightarrow(1)$, let $W$ be a standard $\tilde{\mathfrak{g}}[\mu]$-module. Then $W$, viewed as $\widetilde{\mathfrak{g}^{[0]}}$-module, is a direct sum of standard $\widetilde{\mathfrak{g}^{[0]}}$-module. Since the restriction of $Y^{\sigma}_W(\cdot, x)$ to $L_{\mathfrak{g}^{[0]}}(0, \ell)$ gives $W$ an $\exp(\mbox{ad}\; e)$-twisted $L_{\mathfrak{g}^{[0]}}(0, \ell)$-module structure, if $e_{\theta} \in \mathfrak{g}^{[0]}$, by Theorem \ref{mainthmpaper} $(1)$, $Y^{\sigma}_0(e_{\theta}(-1)^{\ell + 1}, x) = 0$ on $W$. If $e_{\theta} \notin \mathfrak{g}^{[0]}$, then by Lemma \ref{exception}, $Y^{\sigma}_0(e_{\theta}(-1)^{\ell + 1}, x)$ acts on $W$ as $Y^{\mu}(e_{\theta}(-1)^{\ell + 1}, x)$ (see (\ref{affineiso})), thus is $0$.

We only need to prove $(1)\Rightarrow(3)$, let $W$ be a $\sigma$-twisted $L_{\mathfrak{g}}(0, \ell)$-module. Then $W$ can be viewed as an $\exp(\mbox{ad}\; e)$-twisted $L_{\mathfrak{g}^{[0]}}(0, \ell)$-module. By Theorem \ref{mainthmpaper} $(1)$ and $(4)$, $W$ is a direct sum of standard $\widetilde{\mathfrak{g}^{[0]}}$-module. Thus the root vectors of $\mathfrak{g}$ contained in $\mathfrak{g}^{[0]}$ act nilpotent on $W$. Let $f$ be a root vector in $\mathfrak{g}^{[i]}$ for $i \neq 0$. Since $Y^{\mu}(f(-1){\bf 1}, x)$ acts on $W$ as $Y^{\sigma}_0(f(-1){\bf 1}, x)$, by Lemma \ref{exception}, $Y^{\mu}(f(-1){\bf 1}, x)^{\ell + 1}$ is $0$ on $W$. Therefore, $W$ is a standard $\tilde{\mathfrak{g}}[\mu]$-module. \epfv

As a consequence, we prove the $\sigma$-twisted twisted Zhu's algebra for $L_{\mathfrak{g}}(0, \ell)$ is independent of the inner automorphism involved in $\sigma$.

\begin{cor}
The associative algebra $A_{\sigma}(L_{\mathfrak{g}}(0, \ell))$ is semisimple and $A_{\sigma}(L_{\mathfrak{g}}(0, \ell)) \simeq A_{\mu}(L_{\mathfrak{g}}(0, \ell))$.
\end{cor}
\pf From Theorem \ref{mainthmpaper} $(4)$, $A_{\sigma}(L_{\mathfrak{g}}(0, \ell))$ viewed as a left module for $A_{\sigma}(L_{\mathfrak{g}}(0, \ell))$ is completely reducible. Thus $A_{\sigma}(L_{\mathfrak{g}}(0, \ell))$ is a semisimple algebra. The associative algebra isomorphism $A_{\sigma}(L_{\mathfrak{g}}(0, \ell)) \simeq A_{\mu}(L_{\mathfrak{g}}(0, \ell))$ follows from Wedderburn's Theorem and Theorem \ref{ind}. \epfv

We conjecture that it is generally true that for a vertex operator algebra and a general automorphism $g$, the $g$-twisted Zhu's algebra is independent of the inner automorphism involved in $g$:
\begin{conj}
Let $V$ be a vertex operator algebra, let $g$ be an automorphism of $V$ and $h$ be an inner automorphism of $V$. Then
\begin{itemize}
\item[(1)]There is a 1-1 correspondence between the simple objects in the category of lower bounded logarithmic $g$-twisted modules and the simple objects in the category of lower bounded logarithmic $gh$-twisted modules.
\item[(2)]The $g$-twisted Zhu's algebra is isomorphic to the $gh$-twisted Zhu's algebra.
\end{itemize}
\end{conj}

\noindent {\small \sc Department of Mathematics, University of Notre Dame,
278 Hurley Building, Notre Dame, IN 46556}

\noindent {\em E-mail address}:  {\tt  jyang7@nd.edu}

\end{document}